\documentclass[10pt]{amsart}
\usepackage{latexsym}
\usepackage{amssymb}
\usepackage{amsmath}
\usepackage{hyperref}
\usepackage{mathrsfs}

\begin{document}

\title[Saunders Mac Lane]
{Saunders Mac Lane,\break the Knight of Mathematics}
\author{S. S. Kutateladze}
\begin{abstract}
This is a short obituary of Saunders Mac~Lane
(1909--2005).
\end{abstract}
\date{June 11, 2005}


\address[]{
Sobolev Institute of Mathematics\newline
\indent 4 Koptyug Avenue\newline
\indent Novosibirsk, 630090\newline
\indent RUSSIA
}
\email{
sskut@member.ams.org
}
\maketitle

San Francisco and April 14, 2005  form the
terminal place and date of the  marvellous almost centennial life
of the prominent American mathematician
Saunders Mac Lane who shared with
Samuel Eilenberg (1913--1998) the honor of creation of
category theory which ranks among the most brilliant,
controversial, ambitious, and heroic  mathematical achievements
of the 20th century.

Category theory, alongside set theory, serves as
a universal language of modern mathematics.
Categories, functors, and natural transformations
are  widely used in all areas of mathematics, allowing us
to look uniformly and consistently on various
constructions and formulate the general properties
of diverse structures.
The impact of category theory is irreducible to the
narrow frameworks of its great expressive conveniences.
This theory has drastically changed
our general outlook on the foundations of mathematics
and widened the room of free thinking in mathematics.

Set theory, a~great and ingenious creation of Georg Cantor,
occupies in the common opinion of the 20th century the place of
the sole solid base  of  modern  mathematics.
Mathematics becomes  sinking into
a section of the Cantorian set theory.
Most active mathematicians, teachers, and philosophers
consider as obvious and undisputable the thesis
that mathematics cannot be grounded on anything but set theory.
The set-theoretic stance transforms  paradoxically
into an ironclad dogma, a clear-cut forbiddance of
thinking (as  L.~Feuerbach once put it wittily).
Such an indoctrinated view of the foundations of mathematics
is  false and  conspicuously contradicts the leitmotif, nature, and
pathos of the essence of all creative contribution of Cantor
who wrote as far back as in 1883
that  ``denn das {\it Wesen\/} der {\it Mathematik\/} liegt grerade
in ihrer{\it  Freiheit}.''

It is category theory that one of the most ambitious projects
of the 20th century mathematics was realized within in the 1960s,
the project of socializing set theory.
This led to topos theory  providing  a profusion of
categories of which classical set theory is an ordinary
member. Mathematics has acquired infinitely many new
degrees of  freedom. All these  rest on
category theory originated with
the article by
Mac Lane and Eilenberg ``General Theory of Natural Equivalences,''
which was presented to the American Mathematical Society on
September 8, 1942 and published in 1945 in the {\it Transactions of the AMS}.

Mac Lane authored or coauthored more than 100 research papers and
6~books:
{\scshape A~Survey of Modern Algebra} (1941, 1997; with Garrett Birkhoff);
{\scshape Homology} (1963); {\scshape Algebra} (1967; with Garrett Birkhoff); {\scshape Categories for
the Working Mathematician} (1971, 1998); {\scshape Mathematics, Form and Function}
(1985);
{\scshape Sheaves in Geometry and Logic: A First Introduction to Topos Theory}
(1992; with Ieke  Moerdijk).

Mac Lane was the advisor of 39 Ph.D. theses.
Alfred Putman, John Thompson, Irving Kaplansky, Robert Solovay, and many other distinguished scientists
are listed as his students. He was elected to the National Academy of Sciences
of the USA  in 1949 and received  the National Medal of Science, the highest
scientific award of the USA in 1989.
Mac Lane served as vice-president of the National Academy of Sciences
and the American Philosophical Society. He was elected as president of
the American Mathematical Society and Mathematical Association of America.
He contributed greatly to modernization of the teaching programs in mathematics.
Mac Lane received many signs of honor from the leading universities
of the world and possessed an impressive collection of mathematical awards and prizes.
Mac Lane became a living legend of the science of the USA.

Mac Lane was born  on August  4, 1909 in Norwich near Taftville,
Connecticut in the family of a~Congregationalist minister
and was christened as Leslie Saunders MacLane.
The name Leslie was suggested by his nurse, but his mother disliked the name.
A month later,  his father put a hand on the head of the son, looked up to the
God, and said: ``Leslie forget.''
His father and uncles changed the spelling of their surname and began to
write MacLane instead of  MacLean in order to avoid sounding Irish.
The space in Mac Lane was added by Saunders himself at request of
his first wife Dorothy. That is how Mac Lane narrated about his name
in {\scshape A Mathematical Biography\/} which was published soon after his death.

Saunders's father  passed away when the boy was 15 and
it was  Uncle John who supported the boy and paid for
his education in  Yale. Saunders was   firstly fond of chemistry
but everything changed after acquaintance with differential
and integral calculus by the textbook of
Longley and Wilson (which reminds of the later book by
Granville, Smith, and Longley). The university years revealed
Mac Lane's attraction to philosophy and foundations of mathematics.
He was greatly impressed by the brand-new three volumes
by Whitehead and Russell, the celebrated {\scshape Principia Mathematica}.
The mathematical tastes of Mac Lane were strongly influenced by
the lectures of a young assistant professor
Oystein Ore, a Norwegian mathematician from the Emmy Noether's school.
After graduation from Yale, Mac Lane continued education in the
University of Chicago. At that time he was very much influenced by
the personalities and research of
Eliakim Moore, Leonard Dickson, Gilbert Bliss, Edmund Landau,
Marston Morse, and many others.
Mac Lane was inclined to wrote a~Ph.D. thesis in
logic but this was impossible in Chicago and so
Saunders decided to continue education in G\"ottingen.

The stay in Germany in 1931--1933 was decisive for the maturity of
Mac Lane's gift and personality. Although David Hilbert had retired,
he still delivered weekly lectures on
philosophy and relevant general issues.
The successor of Hilbert was Hermann Weyl who had recently arrived
from Z\"urich and was in the prime of his years and talents.
Weyl advised  Saunders to attend the lectures on
linear associative algebras by~Emmy Noether
whom Weyl called ``the equal of each of us.''
In the Mathematical Institute  Mac Lane
met and boiled with Edmund Landau, Richard Courant,
Gustav Herglotz, Otto Neugebauer, Oswald Teichm\"uller, and many others.
Paul Bernays became the advisor of Mac Lane's Ph.D thesis
 ``Abbreviated Proofs in Logic Calculus.''

The Nazis gained power in Germany in February 1933.
The feast of antisemitism started immediately and
one of the first and fiercest strokes fell upon the
Mathematical Institute.  The young persons are welcome to read
as an antidote Mac Lane's masterpiece 

``Mathematics at G\"ottingen under the Nazis''

in~the {\it Notices of the AMS}, {\bf 42}:10, 1134--1138 (1995).

In the fall of 1933  Mac Lane returned to the States  with Dorothy Jones  Mac Lane
whom he had married recently in Germany.
The further academic career of Mac Lane was mainly tied with Harvard
and since 1947 with Chicago.

To evaluate the contribution of  Mac Lane to mathematics is an easy and pleasant task.
It suffices to cite the words  A.~G. Kurosh, a renowned Russian professor
of Lomonosov State University.
In the translator's preface to the Russian edition of
the classical {\scshape Homology\/} book, Kurosh wrote:

\begin{itemize}
\item[]{\it\small\quad
The author of this book, a professor of Chicago University,
is one of the most prominent American algebraists and topologists.
His role in homological algebra as well as category theory
is the role of one of the founders of this area.}
\end{itemize}

Homological algebra implements a marvelous
project of algebraization of  topological spaces
by assigning to such a space $X$  the sequence
of (abelian) homology groups $H_n(X)$. Moreover, each
continuous map
$f: X\to Y$  from $X$ to $Y$ induces a family of
homomorphisms of the homology groups $f_n: H_n(X)\to H_n(Y)$.
The aim of homological algebra consists in calculation of
homologies.

In his research into homological algebra and category theory
Mac Lane cooperated with~Eilenberg whom he met in~1940.
Eilenberg had arrived from Poland two years earlier.
He saw the affinity of the algebraic calculations of
Mac Lane with those he encountered in algebraic topology.
Eilenberg offered cooperation to Mac Lane. The union of
Eilenberg and Mac Lane lasted for 14 years and resulted in
15 joint papers which noticeably changed the mathematical appearance of
the 20th century.

The pearl of this cooperation was category theory.
Mac Lane always considered  category theory ``a~natural
and perhaps inevitable aspect of the 20th
century mathematical emphasis on axiomatic and abstract
methods---especially as those methods when involved in
abstract algebra and functional analysis.''
He stressed that even if  Eilenberg and he did not
propose this theory it will necessarily appear in the works of other
mathematicians.  Among these potential inventors of the new conceptions
Mac Lane listed Claude Chevalley, Heinz Hopf, Norman Steenrod, Henri Cartan,
Charles Ehresmann, and John von Neumann.

In  Mac Lane's opinion, the conceptions of category theory
were close to the methodological principles of the
project of Nicholas Bourbaki.  Mac Lane was sympathetic with the project
and was very close to joining in but this never happened
(the main obstacles were in linguistic facilities). However,
even the later membership of Eilenberg in the Bourbaki group
could not overcome a shade of slight disinclination and repulsion.
It turned out impossible to ``categorize Bourbaki''
with a theory of non-French origin as Mac Lane had once phrased the matter
shrewdly and elegantly.
It is worth noting in this respect that the term
``category theory'' had roots in the mutual interest of its authors in
philosophy and, in particular, in the works of Immanuel Kant.

Set theory rules in the present-day mathematics.
The buffoon's role of ``abstract nonsense''
is assigned in mathematics to category theory.
History and literature demonstrate to us that
the relations between the ruler and the jester
may be totally intricate and unpredictable.
Something very similar transpires in  the interrelations
of set theory and category theory and the dependency
of one of them on the other.

From a logic standpoint, set theory and category theory
are instances of a first order theory. The former
deals with sets and the membership relation between them.
The latter speaks of objects and morphisms (or arrows).
Of course,  there is no principle difference between the atomic formulas
$a\in b$ and $a\to b$. However, the precipice in meaning
is abysmal between the two concepts that are formalized by
the two atomic formulas. The stationary universe of
Zermelo-Fraenkel, cluttered up  with uncountably many copies
of equipollent sets confronts the free world of categories,
ensembles of arbitrary nature that are determined by the
dynamics of their transformations.

The individual dualities of set theory, dependent on
the choice of particular realizations of the pairs of objects under study,
give up their places to the universal  {\it natural transformations\/}
of category theory.
One of the most brilliant achievements of category theory
was the development of axiomatic homology theory. Instead of the
homological diversity for topological spaces
(the simplicial homology for a polyhedron,
singular and \v Cech homology, Vietoris homology, etc.)
Eilenberg and Steenrod suggested as far back as in 1952
the new understanding of each homology or cohomology theory
as a functor from the category of spaces under consideration
to the category of groups.
The axiomatic approach to defining such a functor
radically changed the manner of further progress in homological algebra
and algebraic topology. The study of the homology of
Eilenberg--Mac Lane spaces and the method of acyclic models
demonstrated the strength of the ideas of category theory
and led to universal use of simplicial sets in
$K$-theory and sheaves.

In 1948 Mac Lane proposed the concept of abelian category
abstracting the categories of abelian groups and vector spaces which
played key roles in the first papers on axiomatic homology theory.
The abelian categories were rediscovered in 1953
and became a major tool in research into homological algebra
by Cartan, Eilenberg, and their followers.

Outstanding advances in category theory
are connected with the names of
Alexander Grothendieck and F.~William Lawvere. Topos theory,
their aesthetic creation, appeared in the course of ``point elimination''
called upon by the challenge of invariance of the
objects we study in mathematics.
It is on this road that we met the conception of variable sets
which led to the notion of topos and
the understanding of the social medium of set-theoretic models.

A category is called an {\it elementary topos\/} provided that
it is cartesian closed and has a~suboject classifier.
The sources of toposes lie in the theory of sheaves and Grothendieck topology.
Further progress of the concept of topos
is due to search for some category-theoretic
axiomatization of set theory
as well as study into forcing and the nonstandard set-theoretic models
of Dana~Scott, Robert Solovay, and Petr Vop\v enka.
The new frameworks provide a natural place
for the Boolean valued models that are viewed now
the toposes with Aristotle logic which pave king's ways
to the solution of the  problem of the continuum by
Kurt G\"odel and Paul Cohen. These toposes are now
the main arena of Boolean valued analysis.

Bidding farewell to Mac Lane, reading his sincere and openhearted autobiography,
enjoying his vehement polemics with Freeman J.~Dyson, and perusing his
deep last articles on general mathematics, anyone cannot help but
share his juvenile  devotion and love of mathematics and its creators.
His brilliant essays
 ``Despite Physicists, Proof Is Essential in Mathematics''
and ``Proof, Truth, and Confusion''
form an anthem of mathematics which is only possible by proof.

\begin{itemize}
\item[]{\it\small\quad
 Let me  summarize where we have come.  As with any branch of learning, the
 real substance of mathematics resides in the ideas.  The ideas  of
 mathematics are those which can be formalized and which have  been
 developed  to  fit issues arising in science or in human activity.
 Truth in mathematics is approached by way of proof  in  formalized
 systems.   However,  because  of  the  paradoxical  kinds of
 self-reference exhibited by the barn door and Kurt G\"odel,  there  can
 be no  single formal system which subsumes all mathematical proof.
 To  boot,  the  older  dogmas  that  ``everything  is  logic''   or
 ``everything is a set''  now  have  competition---``everything  is a
 function.''  However, such questions of foundation are  but  a very
 small part of mathematical activity,  which  continues  to  try to
 combine the right ideas to attack substantive problems.   Of these
 I have touched on only a few examples:  Finding all simple groups,
 putting groups together by extension, and  characterizing  spheres
 by their connectivity.  In such cases,  subtle  ideas,  fitted  by
 hand to the problem, can lead to triumph.
\item[]\quad
 Numerical and mathematical methods can be used for  practical
 problems.  However, because of political pressures, the desire for
 compromise, or the simple  desire  for  more  publication,  formal
 ideas may be applied in practical cases where the ideas simply  do
 not  fit.   Then  confusion  arises---whether  from   misleading
 formulation  of  questions  in  opinion  surveys,   from  nebulous
 calculations  of  airy  benefits, by regression, by extrapolation,
 or otherwise.  As the case of fuzzy sets indicates, such confusion
 is  not  fundamentally  a  trouble  caused  by  the  organizations
 issuing reports, but is occasioned by academicians making careless
 use of good ideas where they do not fit.
\item[]\quad
 As Francis Bacon once said, ``Truth ariseth more readily from
 error  than  from  confusion.''   There  remains  to  us, then, the
 pursuit  of  truth,  by  way  of proof, the concatenation of those
 ideas which fit, and the beauty which results when they do~fit.}
\end{itemize}

So wrote   Saunders Mac Lane,
a great genius, creator, master, and servant of mathematics.
His unswerving devotion to the ideals of truth and free thinking
of our ancient science made him the eternal and tragicomical
mathematical Knight of the Sorrowful Figure...

\end{document}